\pdfoutput=1
\documentclass[3pt]{elspreprint}

\usepackage{amssymb}

\usepackage{times}
\usepackage{mathptmx}

\usepackage{listings}
\usepackage{color}
\definecolor{dkgreen}{rgb}{0,0.6,0}
\definecolor{gray}{rgb}{0.5,0.5,0.5}
\definecolor{mauve}{rgb}{0.58,0,0.82}

\journal{Computer Methods in Applied Mechanics and Engineering }

\begin{document}

\begin{frontmatter}

\title{A simple approach to the numerical simulation with trimmed CAD surfaces}

\author[label1,label2]{G. Beer}
\author[label1]{B. Marussig and J.Zechner}
\address[label1]{Institute for Structural Analysis, Graz University of Technology, Lessingstrasse 25, Graz, Austria, e-mail: gernot.beer@tugraz.at}
\address[label2]{Centre for Geotechnical and Materials Modelling,
The University of Newcastle, Australia}

\begin{abstract}
In this work a novel method for the analysis with trimmed CAD surfaces
is presented. The method involves an additional mapping step and the
attraction stems from its simplicity and ease of implementation into
existing Finite Element (FEM) or Boundary Element (BEM) software. The method is first verified with classical test examples in structural mechanics. Then two practical applications are presented one using the FEM, the other the BEM, that show the applicability of the method.
\end{abstract}

\begin{keyword}
trimming, CAD, isogeometric analysis.

\end{keyword}

\end{frontmatter}


\section{Introduction}
 
\begin{figure}
\begin{center}
\includegraphics[scale=0.9]{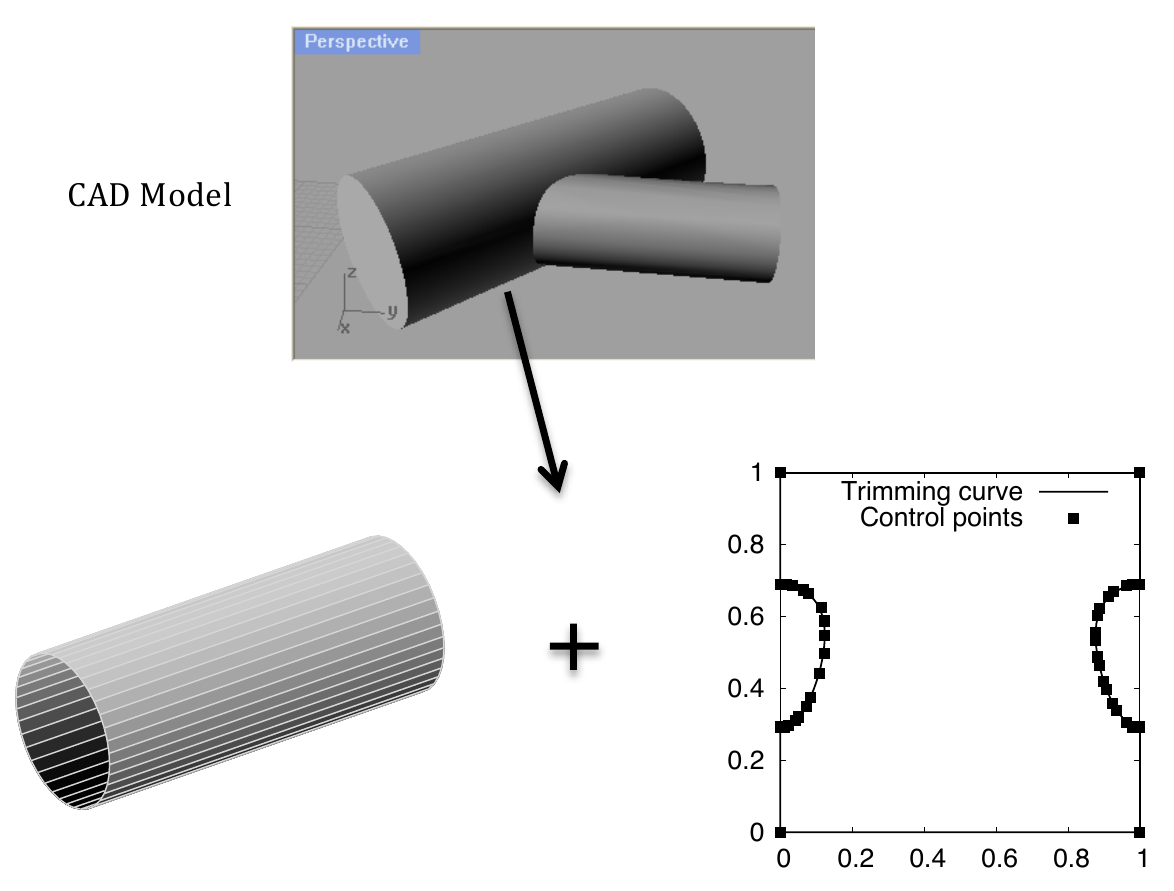}
\caption{Example of IGES output: Tunnel intersection and NURBS surface and trimming information for one tunnel}
\label{IGES}
\end{center}
\end{figure}

The original idea of isogeometric analysis was to use data from CAD programs directly for describing the geometry of the problem without the need to generate a mesh. Geometry data are provided by CAD programs in a standard ASCII format (IGES standard\footnote{available from $http://diyhpl.us/~bryan/papers2/IGES5-3 \_forDownload.pdf$}) and can be read by the analysis program. 
CAD programs describe surfaces using NURBS functions.  For complex surfaces, especially involving a union of surfaces, trimming is used. The data provided in the IGES file contain NURBS information of surfaces (control points, knot vectors and weights), as well as trimming information, if the surface is trimmed. The trimming information is supplied in the global and local coordinate system of the surface to be trimmed.
Figure \ref{IGES} shows an example of an IGES output.

The basic idea of our approach is to use information provided in the IGES file directly for the definition of the geometry and other derived values such as the computation of the Jacobian and first and second derivatives. 
We are aware of only three papers that have addressed the problem of using trimmed CAD data for analysis.
The first two papers (\cite{Kim2009a} and  \cite{Kim2010a}) propose to use a regular grid of elements, that are defined by knot spans. A searching algorithm is employed, that allows to determine how the elements are transected by trimming curves. For the numerical integration trimmed elements are subdivided into triangular subregions. 

The method is very general and can deal with extreme trimming cases such as multiple holes and cases where trimming curves are very close to each other. However, the implementation of the method is not trivial. 
The third paper \cite{Schmidt2012a} deals with the application of trimming to shell surfaces and with the specification of local loading. The method uses a reconstruction of knot spans and control points  and is also generally applicable.

The way our proposed method differs from published ones is that the implementation in existing software is very simple. The application is restricted to cases where a unique mapping exists. In it's current form, this excludes cases with multiple holes that have been analyzed in \cite{Kim2009a} but as will be shown, it is applicable to a number practical problems, especially in the main area of interest of the authors, tunneling. However, expansion of the method is possible.

\section{Trimmed NURBS, double mapping method}

For the explanation of the trimming method, we define two local coordinate systems. One is the local coordinate system of the NURBS surface, used in standard isogeometric analysis, where the local coordinates $u,v$ span from 0 to 1, the other is an additional coordinate system $s,t$, whose coordinates also span form 0 to 1.

\newpage

The new approach to trimming relies on the mapping of an area bordered by 2 trimming curves and straight lines connecting the ends of the curves (defined in the $u,v$ system) into the $s,t$ system. The global coordinates can now be expressed in terms of the local coordinates $s,t$ as:
\begin{equation}
\label{ }
\mathbf{ x}(u,v)= \mathbf{ x}(u(s,t),v(s,t))
\end{equation}

\begin{figure}
\begin{center}
\includegraphics[scale=1.0]{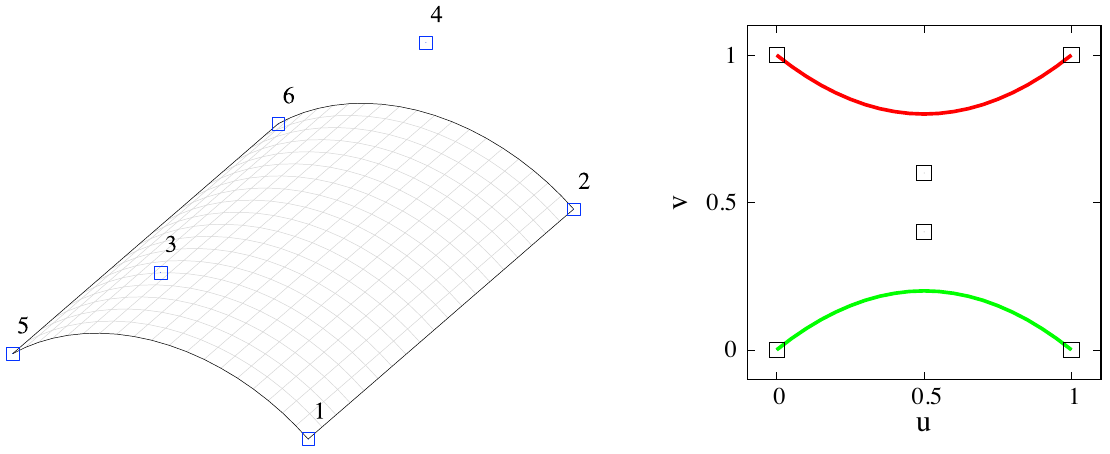}
\caption{Example of trimming. Left: Surface to be trimmed, Right: Two trimming curves with control points}
\label{Triminfo}
\end{center}
\end{figure}

\begin{figure}[h]
\begin{center}
\includegraphics[scale=0.8]{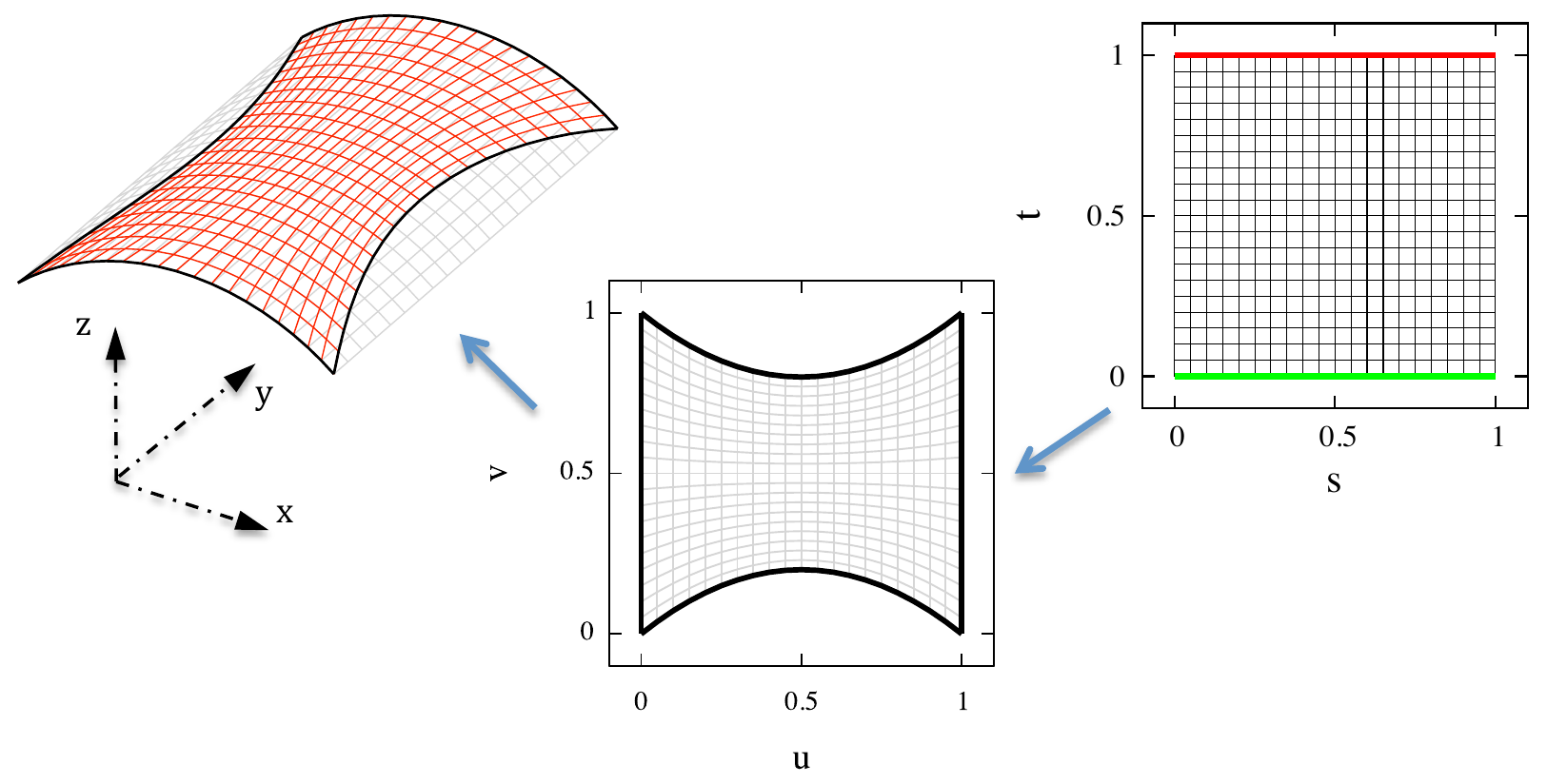}
\caption{Explanation of mapping: right: $s,t$ map, middle: $u,v$ map, left: $x,y,z$ map}
\label{Map}
\end{center}
\end{figure}

This allows all operations, such as computing the derivatives of functions and the numerical integration, to be carried out in the $s,t$ coordinate system. A visual explanation of this procedure is provided in Figures \ref{Triminfo} and \ref{Map}.

We start with mapping the two trimming curves defining the boundaries of the trimmed surface:
\begin{eqnarray}
u^{I}(s)&=&\sum_{n=1}^{N^{I}} R_{n}^{I}(s) \cdot u_{n}^{I} \: ; \: v^{I}(s)=\sum_{n=1}^{N^{I}} R_{n}^{I}(s) \cdot v_{n} ^{I} \\
\nonumber
u^{II}(s)&=&\sum_{n=1}^{N^{II}} R_{n}^{II}(s) \cdot u_{n}^{II} \: ; \: v^{II}(s)=\sum_{n=1}^{N^{II}} R_{n}^{II}(s) \cdot v_{n} ^{II} 
\end{eqnarray}

where $R_{n}^{I}(s)$, $R_{n}^{II}(s)$ are one-dimensional NURBS basis functions defining the trimming curves, $N^{I}$, $N^{II}$ are the number of control points and $u_{n}^{I},v_{n} ^{I}$ and  $u_{n}^{II},v_{n} ^{II}$ are the $u,v$ coordinates of control points. The superscript refers to the curve number ($II$=top; $I$=bottom). As can be seen in Figure \ref{Map} these curves map as straight lines in the $s,t$ coordinate system.

Next we use a linear interpolation between the curves to map the trimmed area:
\begin{eqnarray}
u(s,t)= N_{1}(t)\cdot u^{I}(s) +  N_{2}(t)\cdot u^{II}(s)  \\
v(s,t)= N_{1}(t)\cdot v^{I}(s) +  N_{2}(t)\cdot v^{II}(s)  
\end{eqnarray}

with
\begin{equation}
\label{ }
N_{1}(t)=1 - t \: ; \: N_{2}(t)= t 
\end{equation}

Finally we map from the $u,v$ system to the global $x,y,z$ system by
\begin{equation}
\mathbf{x}=  \sum_{b=1}^{B}\sum_{a=1} ^{A}R_{a,b}^{p,q} (u,v)\mathbf{x}_{a,b}
\label{Gmap}
 \end{equation}

where $R_{a,b}^{p,q} (u,v)$ are NURBS functions describing the untrimmed NURBS surface. $A,B$ is are the number of control points in $u,v$ directions and $\mathbf{x}_{a,b}$ are control point coordinates.

For the computation of the Jacobian matrix we need the first derivatives of \textbf{x} in terms of $s,t$.

The first derivatives of \textbf{x} to $s$ and $t$ are given by
\begin{eqnarray}
 \frac{\partial \mathbf{ x}}{\partial t}=  \frac{\partial \mathbf{ x}}{\partial u} \cdot \frac{\partial u}{\partial t} + \frac{\partial \mathbf{ x}}{\partial v} \cdot \frac{\partial v}{\partial t}\\
\nonumber
 \frac{\partial \mathbf{ x}}{\partial s}=  \frac{\partial \mathbf{ x}}{\partial u} \cdot \frac{\partial u}{\partial s} + \frac{\partial \mathbf{ x}}{\partial v} \cdot \frac{\partial v}{\partial s} 
\end{eqnarray}

where 
\begin{eqnarray}
\frac{\partial u}{\partial s}  & = &  N_{1}(t) \cdot \frac{\partial u^{I}(s)}{\partial s} + N_{2}(t) \cdot \frac{\partial u^{II}(s)}{\partial s} \\
\nonumber
\frac{\partial v}{\partial s}  & = &  N_{1}(t) \cdot \frac{\partial v^{I}(s)}{\partial s} + N_{2}(t) \cdot \frac{\partial v^{II}(s)}{\partial s} \\
\nonumber
\frac{\partial u}{\partial t}  & = & - 1 \cdot u^{I}(s) + 1 \cdot u^{II}(s)\\
\nonumber
\frac{\partial v}{\partial t}  & = & - 1 \cdot v^{I}(s) + 1 \cdot v^{II}(s)
\end{eqnarray}

and
\begin{equation}
\label{ }
\frac{\partial u^{I}(s)}{\partial s}=\sum_{n=1}^{N^{I}} \frac{\partial R_{n}^{I}(s)}{\partial s} \cdot u_{n}^{I} \: ; \: \frac{\partial u^{II}(s)}{\partial s}=\sum_{n=1}^{N^{II}} \frac{\partial R_{n}^{II}(s)}{\partial s} \cdot u_{n} ^{II}
\end{equation}

For Kirchhoff shell analysis we require the second derivatives also, which are given by:

\begin{eqnarray}
 \frac{\partial^{2} \mathbf{ x}}{\partial s^{2}} =  \left( \frac{\partial^{2} \mathbf{ x}}{\partial u^{2}} \cdot \frac{\partial u}{\partial s}  + \frac{\partial^{2} \mathbf{ x}}{\partial u \partial v}  \cdot \frac{\partial v}{\partial s}\right) \cdot \frac{\partial u}{\partial s} 
 +  \frac{\partial \mathbf{ x}}{\partial u} \cdot \frac{\partial^{2} u}{\partial s^{2}} \\ 
 \nonumber
 + \left( \frac{\partial^{2} \mathbf{ x}}{\partial u \partial v} \cdot \frac{\partial u}{\partial s}  + \frac{\partial^{2} \mathbf{ x}}{\partial v^{2}}  \cdot \frac{\partial v}{\partial s}\right) 
 \cdot \frac{\partial v}{\partial s} +  \frac{\partial \mathbf{ x}}{\partial v} \cdot \frac{\partial^{2} v}{\partial s^{2}} \\
 \nonumber
  \frac{\partial^{2} \mathbf{ x}}{\partial t^{2}} =  \left( \frac{\partial^{2} \mathbf{ x}}{\partial u^{2}} \cdot \frac{\partial u}{\partial t}  + \frac{\partial^{2} \mathbf{ x}}{\partial u \partial v}  \cdot \frac{\partial v}{\partial t}\right) \cdot \frac{\partial u}{\partial t} 
 +  \frac{\partial \mathbf{ x}}{\partial u} \cdot \frac{\partial^{2} u}{\partial t^{2}} \\ 
 \nonumber
 + \left( \frac{\partial^{2} \mathbf{ x}}{\partial u \partial v} \cdot \frac{\partial u}{\partial t}  + \frac{\partial^{2} \mathbf{ x}}{\partial v^{2}}  \cdot \frac{\partial v}{\partial t}\right) 
 \cdot \frac{\partial v}{\partial t} +  \frac{\partial \mathbf{ x}}{\partial v} \cdot \frac{\partial^{2} v}{\partial t^{2}} \\
  \nonumber
  \frac{\partial^{2} \mathbf{ x}}{\partial t \partial s} =  \left( \frac{\partial^{2} \mathbf{ x}}{\partial u^{2}} \cdot \frac{\partial u}{\partial s}  + \frac{\partial^{2} \mathbf{ x}}{\partial u \partial v}  \cdot \frac{\partial v}{\partial s}\right) \cdot \frac{\partial u}{\partial t} 
 +  \frac{\partial \mathbf{ x}}{\partial u} \cdot \frac{\partial^{2} u}{\partial t \partial s} \\ 
 \nonumber
 + \left( \frac{\partial^{2} \mathbf{ x}}{\partial u^{2} } \cdot \frac{\partial u}{\partial s}  + \frac{\partial^{2} \mathbf{ x}}{\partial v \partial u}  \cdot \frac{\partial v}{\partial s}\right) 
 \cdot \frac{\partial v}{\partial t} +  \frac{\partial \mathbf{ x}}{\partial v} \cdot \frac{\partial^{2} v}{\partial t \partial s}
\end{eqnarray}

The derivatives of $u$ and $v$ are given by:

\begin{eqnarray}
\frac{\partial^{2} u}{\partial s^{2}}  & = &  N_{1}(t) \cdot \frac{\partial^{2} u_{b}(s)}{\partial s_{2}} + N_{2}(t) \cdot \frac{\partial^{2} u_{t}(s)}{\partial s^{2}} \\
\nonumber
\frac{\partial^{2} v}{\partial s^{2}}  & = &  N_{1}(t) \cdot \frac{\partial^{2} v_{b}(s)}{\partial s^{2}} + N_{2}(t) \cdot \frac{\partial^{2} v_{t}(s)}{\partial s^{2}} \\
\nonumber
\frac{\partial^{2} u}{\partial t^{2}}  & = & 0 \: ; \: \frac{\partial^{2} v}{\partial t^{2}}   =  0 \\
\nonumber
\frac{\partial^{2} u}{\partial s \partial t}  & = & -1\cdot \frac{\partial u_{b}}{\partial s}+1\cdot \frac{\partial u_{t}}{\partial s}\\
\nonumber
\frac{\partial^{2} v}{\partial s \partial t}  & = & -1\cdot \frac{\partial v_{b}}{\partial s}+1\cdot \frac{\partial v_{t}}{\partial s}
\end{eqnarray}

 The method can only be used in conjunction with the concept of geometry independent field approximation, which was first introduced in \cite{BeerBordas2014} and applied to practical problems in   \cite{Beer2013}, \cite{Marussig2013a}, \cite{Marus} \cite{Marussig2014a} and \cite{Xu}. 
 
This means that the description of the geometry is unchanged during refinement and that the field variables are approximated independently. For example, unknown displacements are approximated by
\begin{equation}
\mathbf{u}=  \sum_{b=1}^{B_{u}}\sum_{a=1} ^{A_{u}}N_{a,b}^{p,q} (s,t)\mathbf{u}_{a,b}
\label{Gmap}
 \end{equation}

where $N_{a,b}^{p,q} (s,t)$ are basis functions defined in the $s,t$ coordinate system, $\mathbf{u}_{a,b}$ are parameter values and $A_{u},B_{u}$ are the number of parameters in $s,t$ directions. The basis functions and the integration regions are defined in the local $s,t$ coordinate system and then mapped into the global system. Because the mapping affects the continuity of the functions in the global parameter space, care has to be taken in the choice of the functions and integration regions (see comments later).

\section{Implementation}

The method was implemented into isogeometric FEM and BEM software using Octave and the nurbs toolkit\footnote{available for free download from http://octave.sourceforge.net/nurbs/overview.html}. See \cite{BeerBordas2014} and \cite{Beer2015} for more details about the program.
The main mapping function is shown here:
\begin{lstlisting}
function [uv,vecu]= Map(Knotg,Coefs,tt)
%--------------------------------------
%    Maps points on trimmed surface from s,t to u,v 
%    Input:
%    Knotg   ...  array containing knot vectors of all trimming curves
%    Coefs  ...   array containing control points of all trimming curves
%    tt(1,np),tt(2,np)  ... vector with s,t-coordinates
%    Output:
%    uv(1,np),uv(2,np)      ...   u,v coordintes of points
%    vecu    ...   Vectors in u,v directions 
%-------------------------------------
for ncrv=1:2
%  extract knot values and control points of trimming curves
 [knotu,coefs]= Get_infoc(Knotg,Coefs);
%   mapping of trimming curves from s,t to u,v 
 nurbs= nrbmak(coefs,knotu); dnurbs= nrbderiv(nurbs);
 [uv,vuv]= nrbdeval(nurbs,dnurbs,tt(1,:));
 uvn(:,:,ncrv)= uv(:,:); vuvn(:,:,ncrv)=vuv(:,:);
end 
%  mapping of points on trimmed surface from s,t to u,v
for np=1:columns(tt)
   N1=1-tt(2,np); N2=tt(2,np);
% compute u,v values
  uv(1,np)= N1*uvn(1,np,1) + N2*uvn(1,np,2);
  uv(2,np)= N1*uvn(2,np,1) + N2*uvn(2,np,2);
%  first derivatives
  v1(1)= N1*vuvn(1,np,1) + N2*vuvn(1,np,2);
  v1(2)= N1*vuvn(2,np,1) + N2*vuvn(2,np,2); v1(3)= 0;
  v2(1)= uvn(1,np,2) - uvn(1,np,1);
  v2(2)= uvn(2,np,2) - uvn(2,np,1); v2(3)= 0;
  vecu{1}(:,np)=v1; vecu{2}(:,np)=v2;
end
endfunction;
\end{lstlisting}

To implement the method insert the mapping function before the mapping to global coordinates. For example:
\begin{lstlisting}
[uv,vecu]= Map(Knotr,Coefstr,tt);
[x,vect]= nrbdeval(nurbs,ders,uv);
\end{lstlisting}

\begin{figure}[h]
\begin{center}
\includegraphics[scale=0.4]{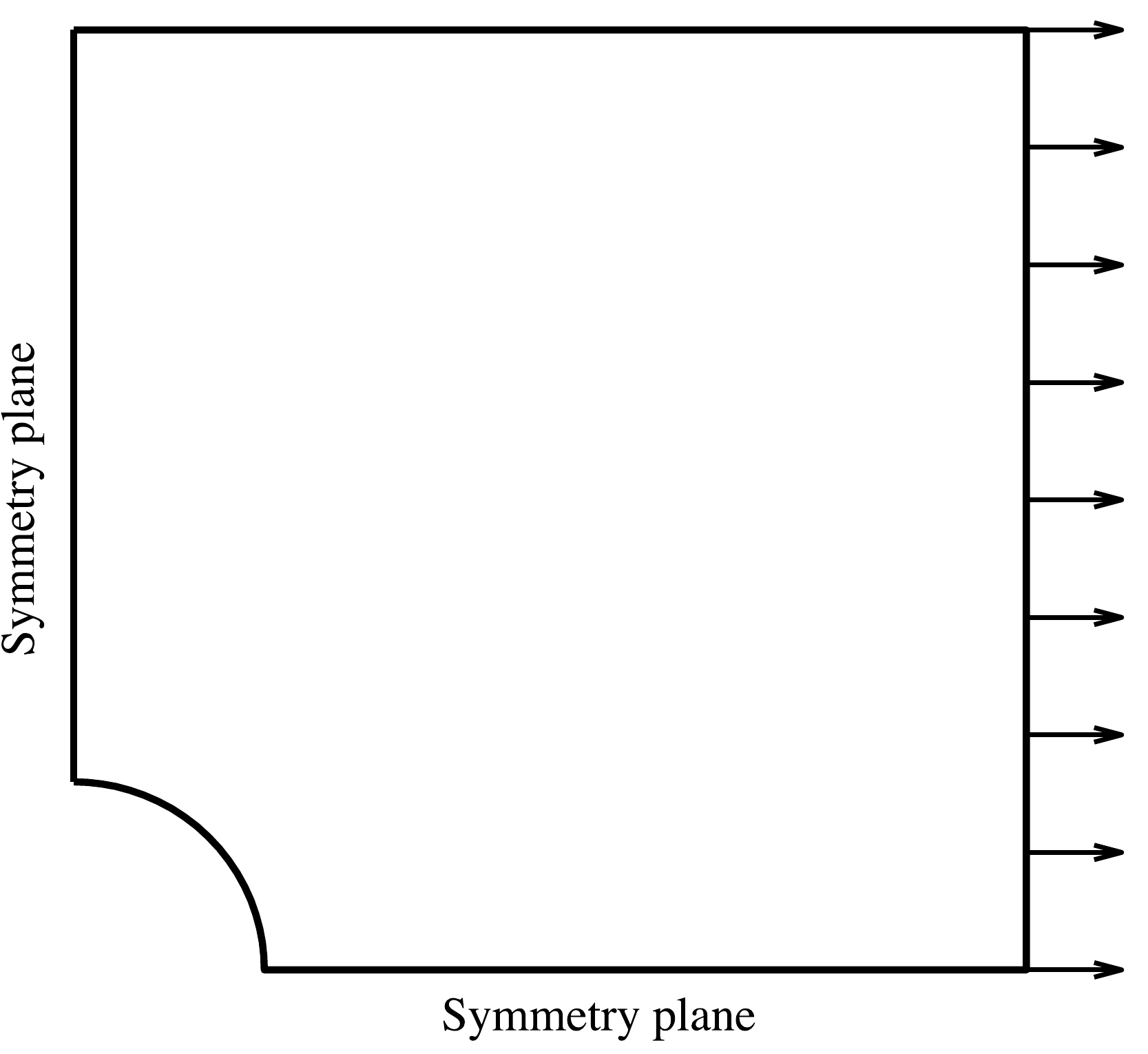}
\caption{Test example 1: Plate with a hole}
\label{Ph}
\end{center}
\end{figure}

\section{Test examples}

We have chosen two test examples, which are the same or similar to the ones used in the cited references on trimming. They relate to the FEM analysis of a plate and a Kirchhoff shell. 
\begin{figure}
\begin{center}
\includegraphics[scale=1.0]{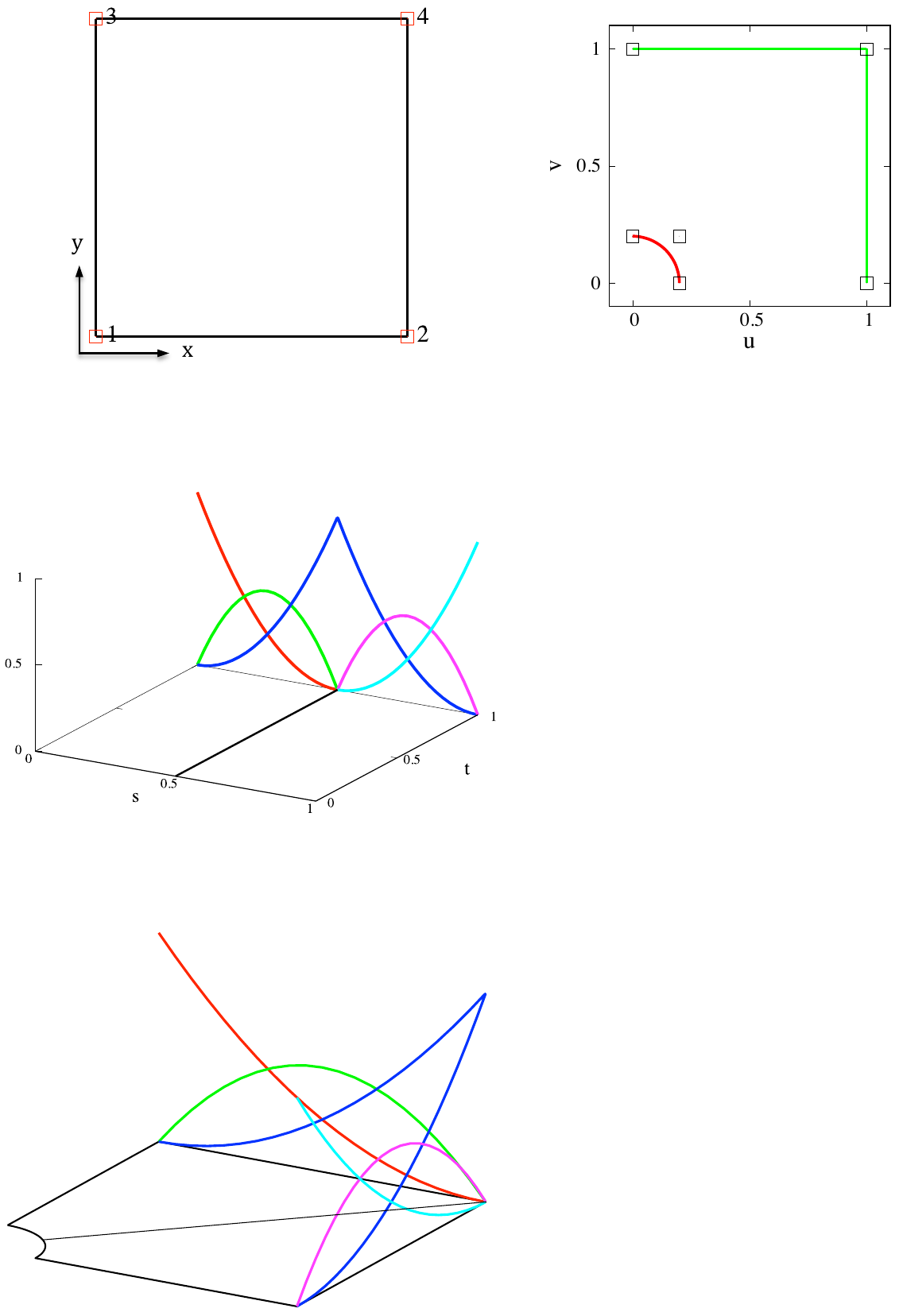}
\caption{Geometry definition of plate with a hole: Top left: untrimmed NURBS surface in global x,y system and Top right: trimming curves in local u,v system. Bottom: trace of basis functions along $t=1$ in the local $s,t$ and in the global $x,y$ coordinate system}
\label{Ph1}
\end{center}
\end{figure}

\subsection{Plate with a hole}

This is the exactly the same test example that has been used by Kim et al \cite{Kim2009a} to verify their trimming method. It relates to a plate with a hole. Exact solutions are available for an infinite plate with a hole. If the correct boundary conditions are applied, then a plate of finite extent can be analyzed and the solution compared with the exact solution. 
Figure \ref{Ph} shows the geometry of the plate with the symmetry planes assumed. The size of the plate  is 5x5 and the radius of the hole is 1. The elastic properties assumed are $E=10^{5}, \nu=0.3$.
The plate is subjected to a tensile horizontal stress of 1 at the right boundary. To simulate an infinite extent, the tractions computed from the exact solution were applied at the top. 
The trimming is explained in Figure \ref{Ph1}. 

The definition of the 2 trimming curves is as follows:
We start with a bilinear NURBS surface describing the plate without the hole and add two trimming curves. The plate is described by:
\begin{lstlisting}
Order: 1
Knot vector: 0 0 1 1
Coefficients (x,y,z,weight):
0 0 0 1
0 1 0 1
1 0 0 1
1 1 0 1
\end{lstlisting}

The trimming curves are defined by:
\begin{lstlisting}
% Curve 1: 
Order: 2
Knot vector: 0 0 0 1 1 1
Coefficients:
0 0.2 0 1
0.2 0.2 0 0.707
0.2 0 0 1

% Curve 2:
Order: 1
Knot vector: 0 0 0.5 1 1
Coeffcients:
0 1 0 1
1 1 0 1
1 0 0 1
\end{lstlisting}

Because one of the trimming curves has only a $C^{0}$ continuity, care has to be taken when choosing the basis functions for the description of the unknown before refinement and the regions used for the integration. In Figure  \ref{Ph1} we show that in this case the integration region is split into two and the continuity of the basis functions is changed to match the one of the mapping function. In this Figure we also show the trace of the basis functions for the approximation of the unknown along $t=1$ in the local $s,t$ coordinate system and in the global $x,y$ coordinate system.
\begin{figure}
\begin{center}
\includegraphics[scale=0.8]{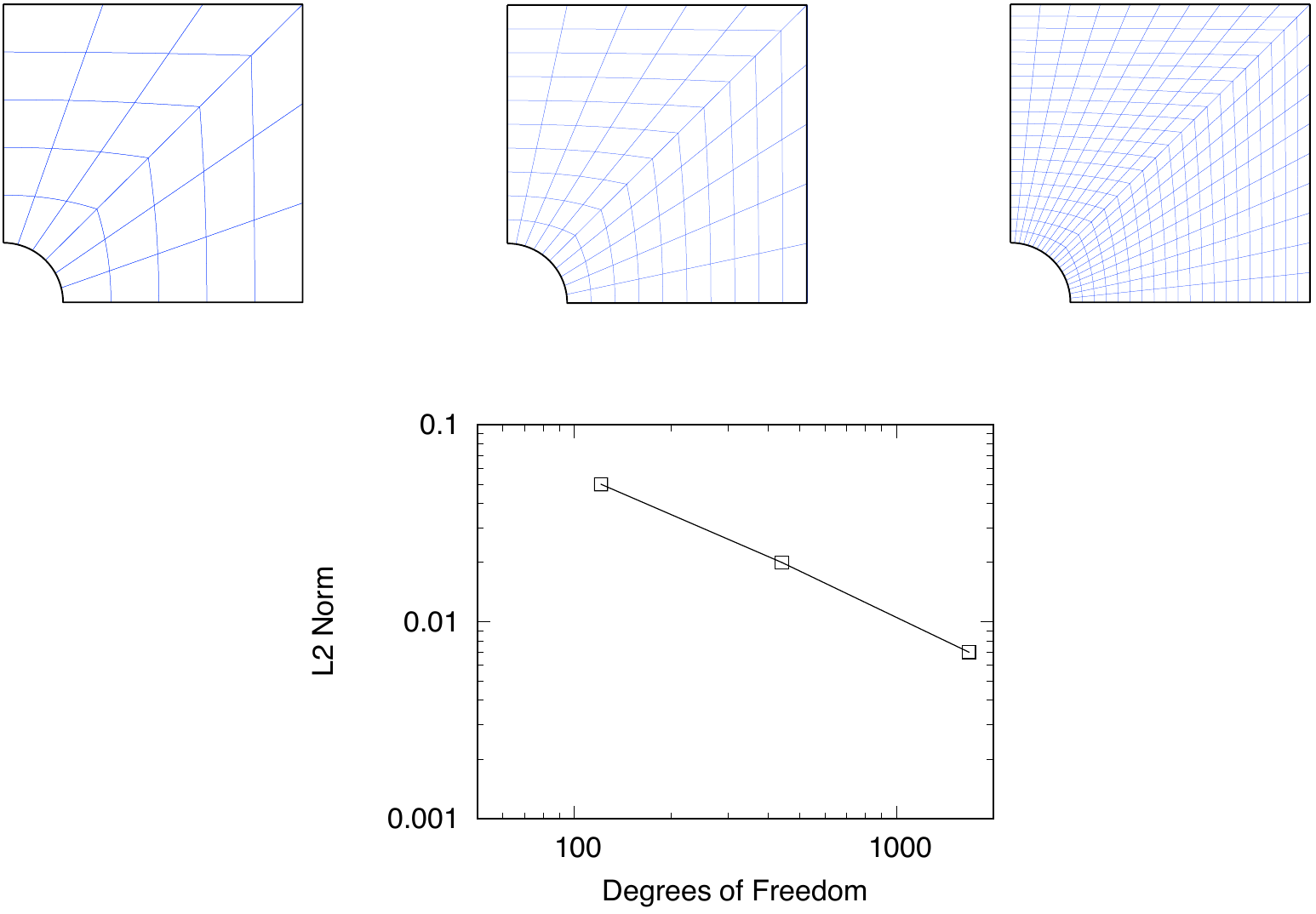}
\caption{The 3 refinement stages investigated and convergence plot}
\label{Ph2}
\end{center}
\end{figure}
\begin{figure}[h]
\begin{center}
\includegraphics[scale=0.45]{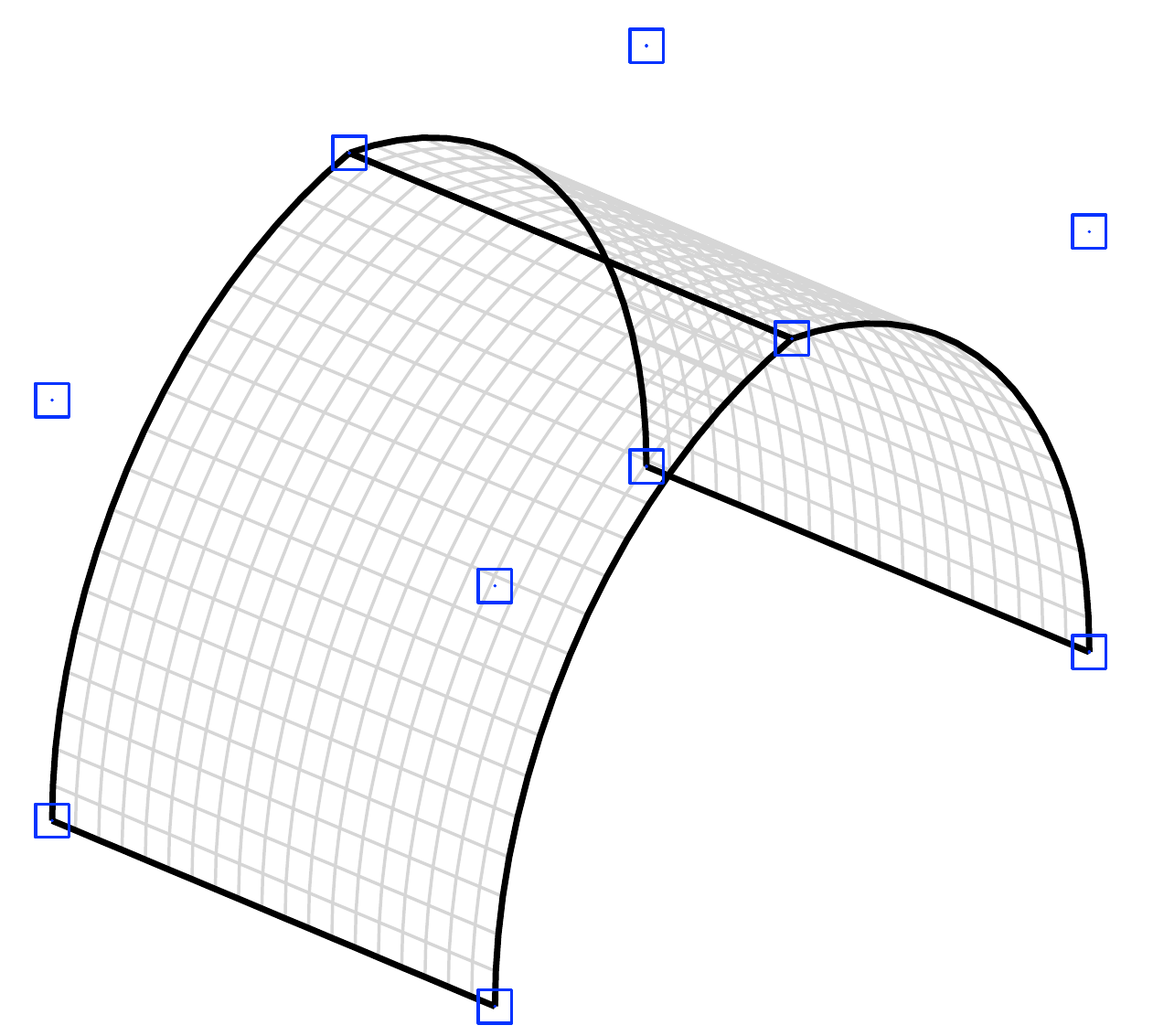}
\caption{Geometry definition of test case with 2 untrimmed NURBS patches, showing control points}
\label{Kh1}
\end{center}
\end{figure}

We analyze three refinements and as expected, the convergence of the $L^{2}$ norm of the stress shown in Figure \ref{Ph2} is similar to the one reported in \cite{Kim2009a} for analyses without trimming.

\subsection{Kirchhoff shell}

This test example is similar to the one used by Schmidt et al in \cite{Schmidt2012a} to test their trimming method. It is the analysis of a Kirchhoff shell. Details of the implementation of isogeometric Kirchhoff shell analysis can be found in \cite{Kiendl}. 
For the test we have implemented our mapping method into software developed as part of a Master Thesis \cite{Fleissner}.

Two analyses are performed: one with two untrimmed and the other with two trimmed NURBS patches.
\begin{figure}
\begin{center}
\includegraphics[scale=1.0]{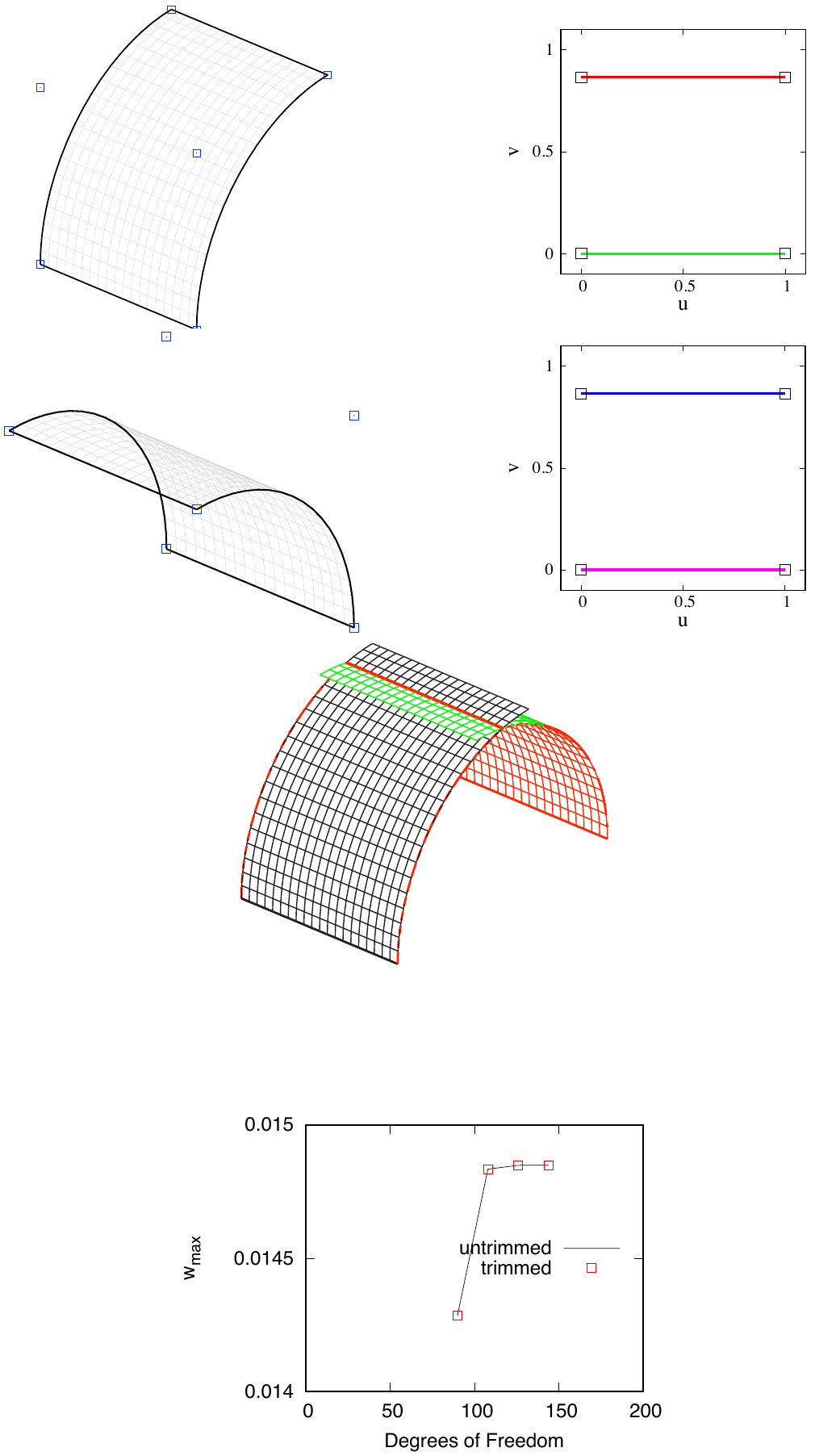}
\caption{Two untrimmed NURBS surfaces with trimming curves, resulting trimmed surface and Convergence of maximum displacement}
\label{Kh2}
\end{center}
\end{figure}

Figure \ref{Kh1} shows the geometry definition with 2 untrimmed NURBS patches. The patches are fixed in three coordinate directions at the base. 

At the top only the translational degrees of freedom are connected (hinged connection). The shell is subjected to selfweight.
Figure \ref{Kh2} shows the case where the two NURBS patches have been extended and then trimmed back to match the geometry of the untrimmed case. Both cases were analyzed and convergence achieved by order elevation. 
The convergence of the maximum deflection for both cases is shown in Figure \ref{Kh2} and as expected, no difference between the cases can be observed.

\section{Practical examples}

Here we show two practical examples that demonstrate that the proposed mapping method, despite its limitations, has practical applications.
One example is in shell theory and is meant to demonstrate how easily the geometry can be changed. The second one relates to the boundary element analysis of a branched tunnel.

\subsection{Example 1: Arched Scordelis roof}

\begin{figure}
\begin{center}
\includegraphics[scale=1.0]{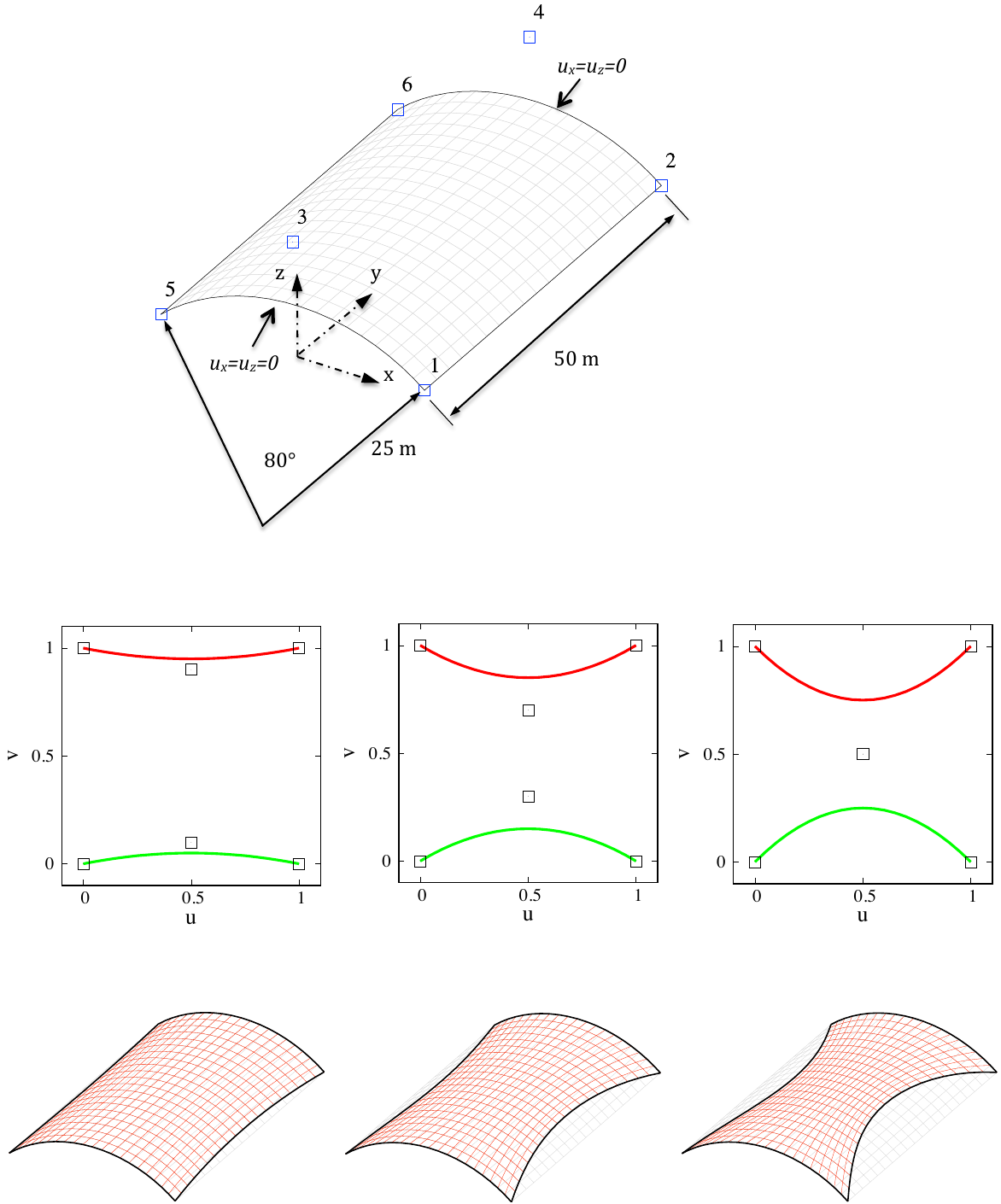}
\caption{Scordelis-Lo roof: Problem definition and the trimming curves for the introduction of different degrees of arching (from left to right geometry 1 to 3)}
\label{Scgeo}
\end{center}
\end{figure}

As a first example we show a variation on a classical example in shell theory, the Scordelis-Lo roof. The geometry of the problem is shown in Figure \ref{Scgeo}.

The shell is subjected to self weight $q$ and the properties are:
\begin{itemize}
  \item $E= 4.32 \cdot 10^8 \: kN/m^2$ , $\nu=0$
  \item $t=0.25$ , $q=90 \: kN/m^2$
\end{itemize}

We change the geometry by introducing some arching and use three different trimming curves with increasing removal of material (Figure \ref{Scgeo}).

\begin{figure}
\begin{center}
\includegraphics[scale=0.65]{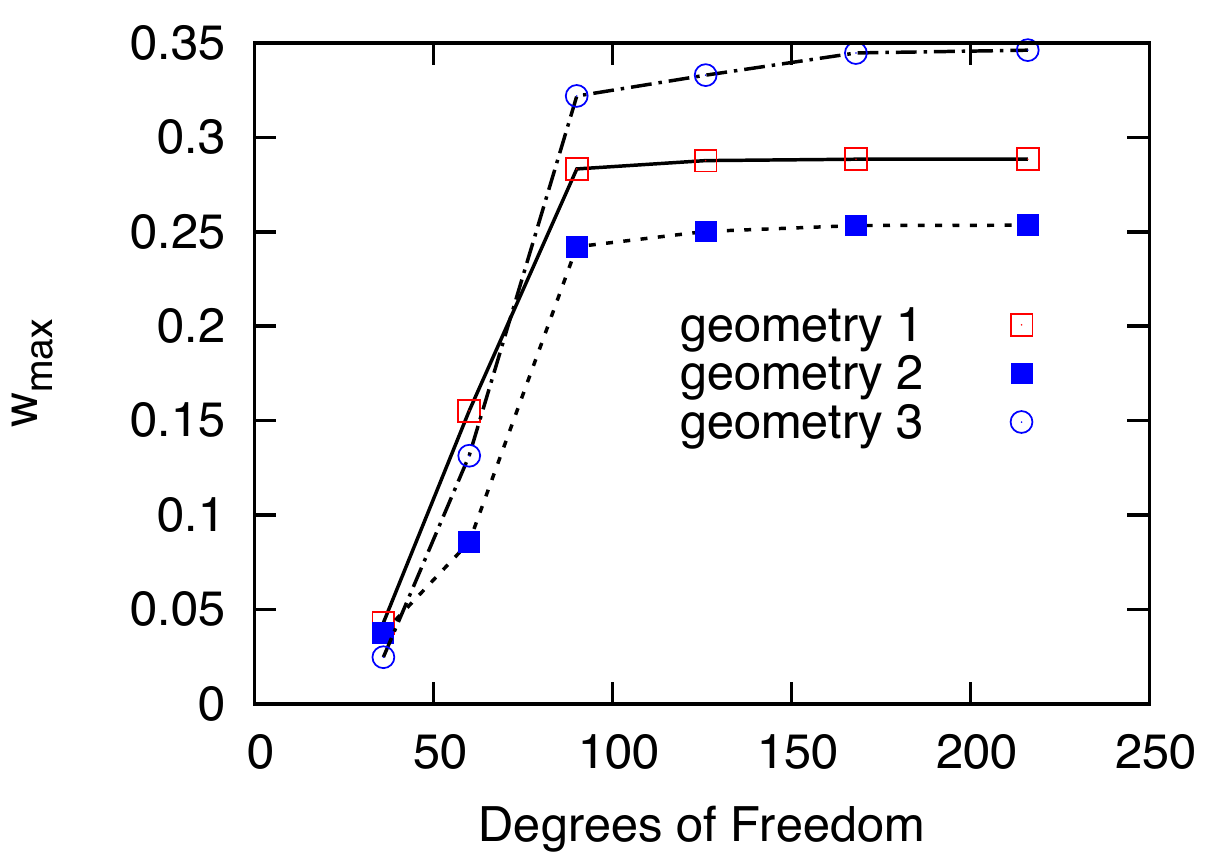}
\caption{Arched Scordelis roof: Convergence of the maximum displacement}
\label{Cvarch}
\end{center}
\end{figure}

\begin{figure}
\begin{center}
\includegraphics[scale=1.0]{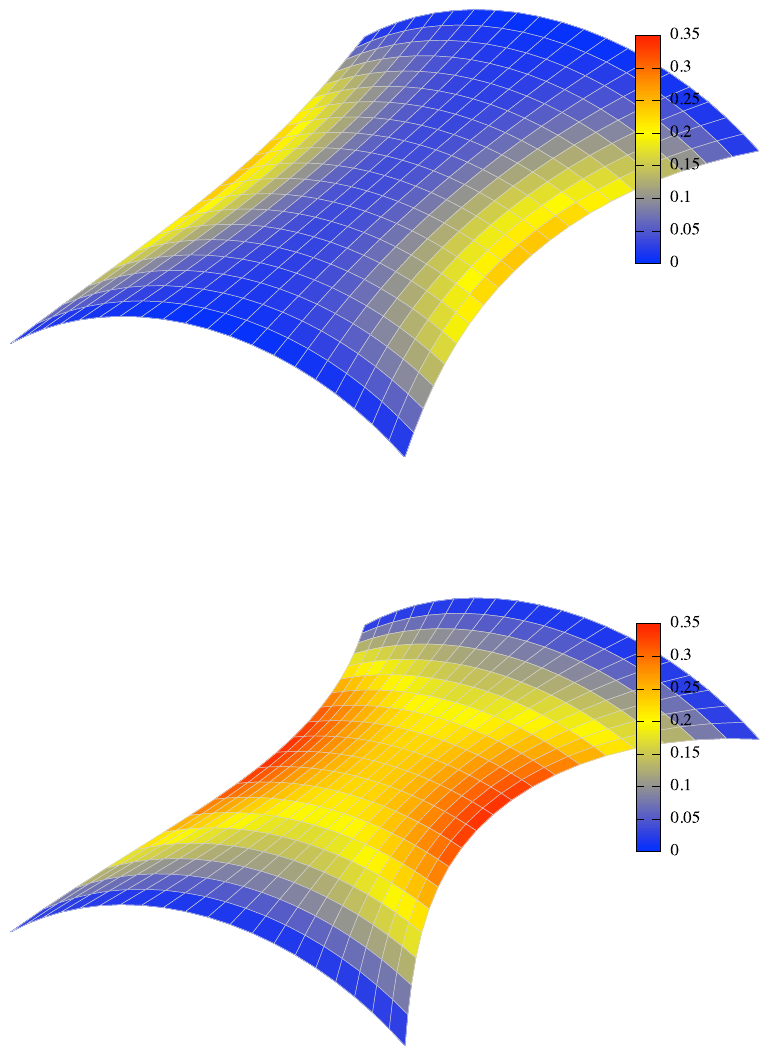}
\caption{Contours of vertical displacement for geometries 2 and 3}
\label{Scardis}
\end{center}
\end{figure}

It can be seen in Figure \ref{Cvarch} that all cases converge. For moderately arched shells the maximum vertical displacement decreases but if more material is removed there is an increase. 

In Figure \ref{Scardis} we can see that for moderate material removal the maximum displacement is localized, which spreads over the shell as more material is removed.

\subsection{Example 2: Tunnel branch}

This example is meant to demonstrate the implementation of the mapping method in a Boundary Element code. The implementation of the isogeometric BEM is discussed in \cite{BeerBordas2014} and  of the mapping method in \cite{Beer2015}. Here only results of the simulation are presented. More details can be found in \cite{Beer2014a}.

\begin{figure}
\begin{center}
\includegraphics[scale=1.0]{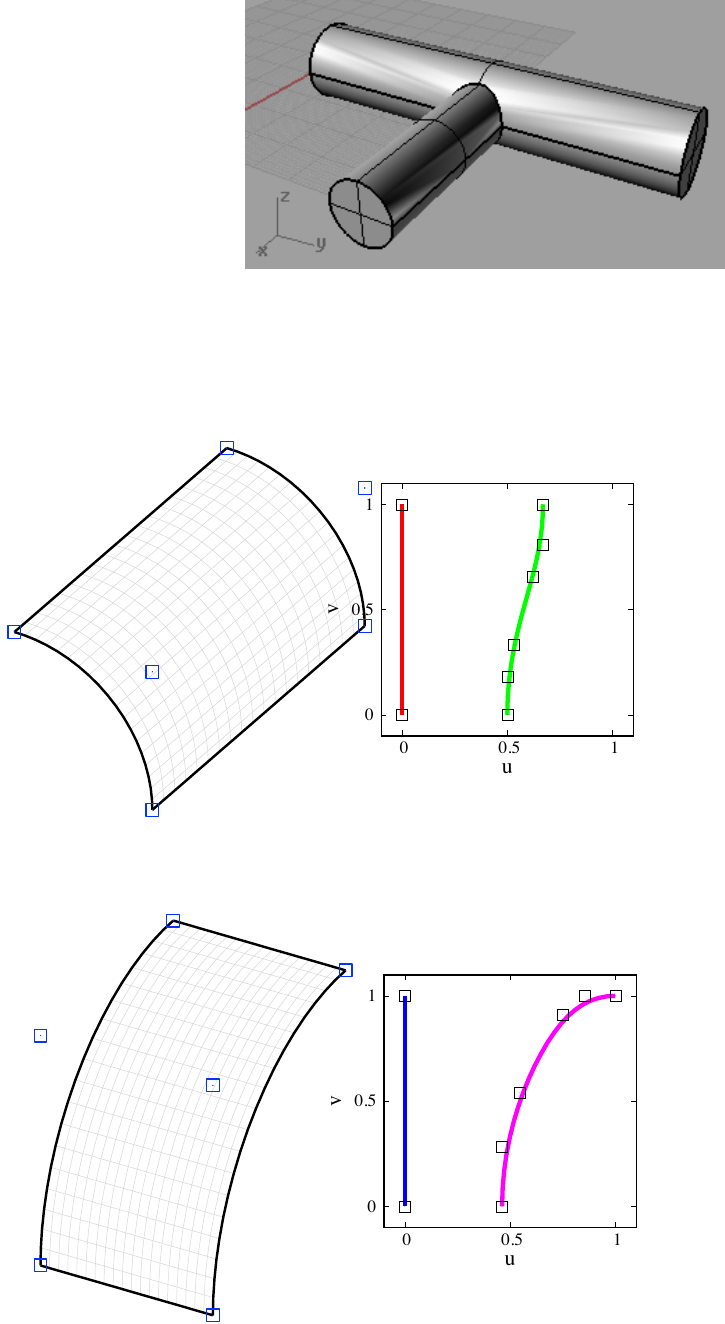}
\caption{CAD model of tunnel intersection and extracted NURBS surface and trimming information}
\label{TB1}
\end{center}
\end{figure}

\begin{figure}
\begin{center}
\includegraphics[scale=0.4]{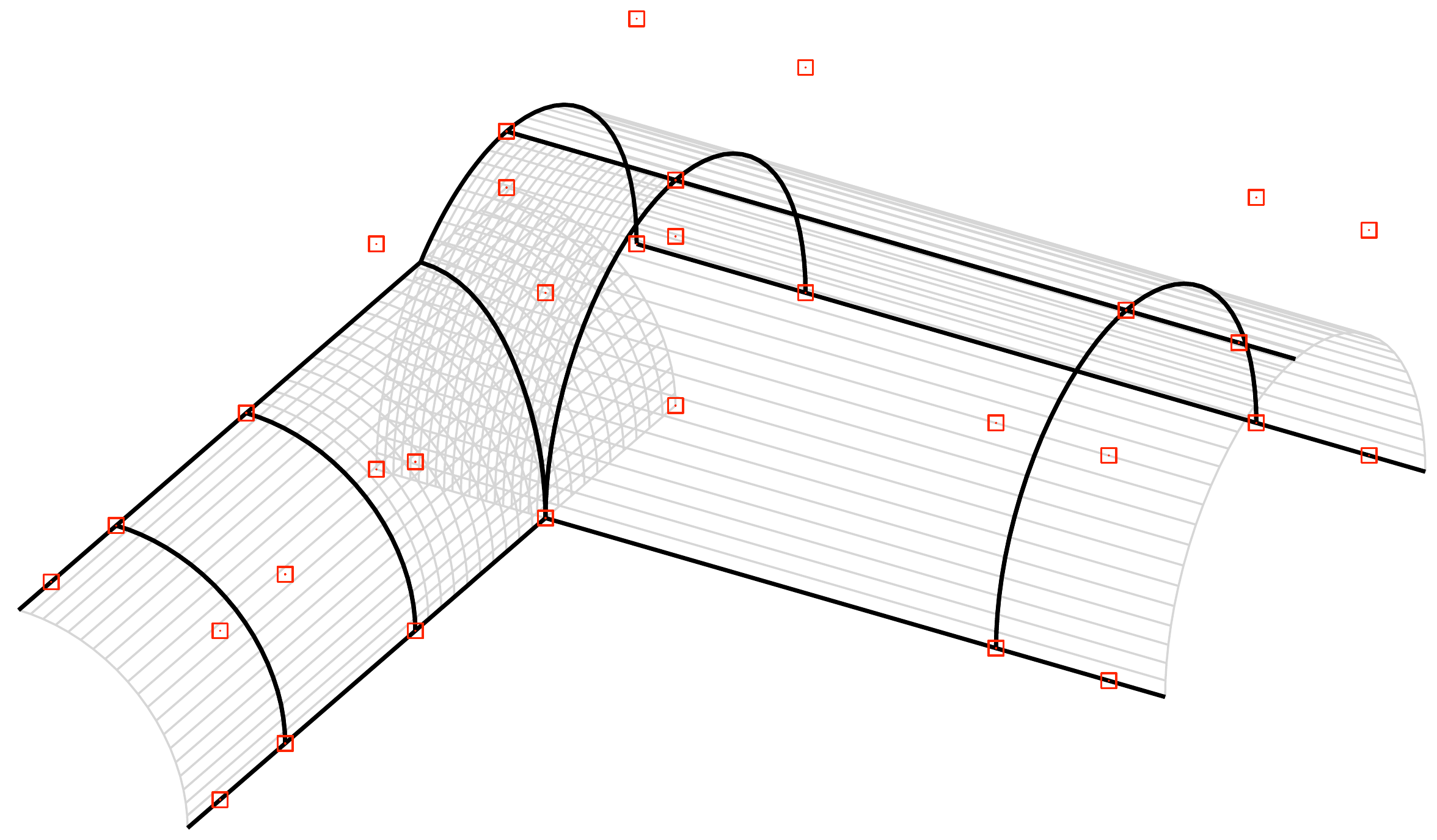}
\caption{Geometry description for the simulation}
\label{TB2}
\end{center}
\end{figure}

 Figure \ref{TB1} shows the Rhino CAD model of a tunnel intersection and the geometrical information extracted from the IGES file produced by the program. This information is used directly for the simulation as follows: The trimmed surfaces are used to model 1/4 of the problem as shown in Figure \ref{TB2}. Two planes of symmetry (about the x-y and x-z planes) are assumed and specially developed infinite plane strain patches (see \cite{BeerBordas2014} and \cite{Beer2015} for details) are used to simulate infinitely long tunnels. The tunnels are assumed to be excavated in one step in an elastic pre-stressed ground. 
 
\

\textbf{Remark}: Since the trimming curves are specified separately for each NURBS patch they will in general not match perfectly and there may be small gaps. The gaps can be reduced to a negligible size by specifying a small tolerance (option "units and tolerances" in Rhino). The other problem is that the parameter spaces of the trimming curves may not match at the interface. The implication of this for the BEM is that the collocation points, whose coordinates are computed using Gervilla abscissa in the local coordinate system of the trimmed NURBS, may not be in the same location. If this issue becomes significant the way of dealing with this would be by using discontinuous collocation (see for example  \cite{Marussig2015a}). For the example presented here the collocation point locations matched within 1 \%, so continuous collocation was used.

\

\begin{figure}
\begin{center}
\includegraphics[scale=0.3]{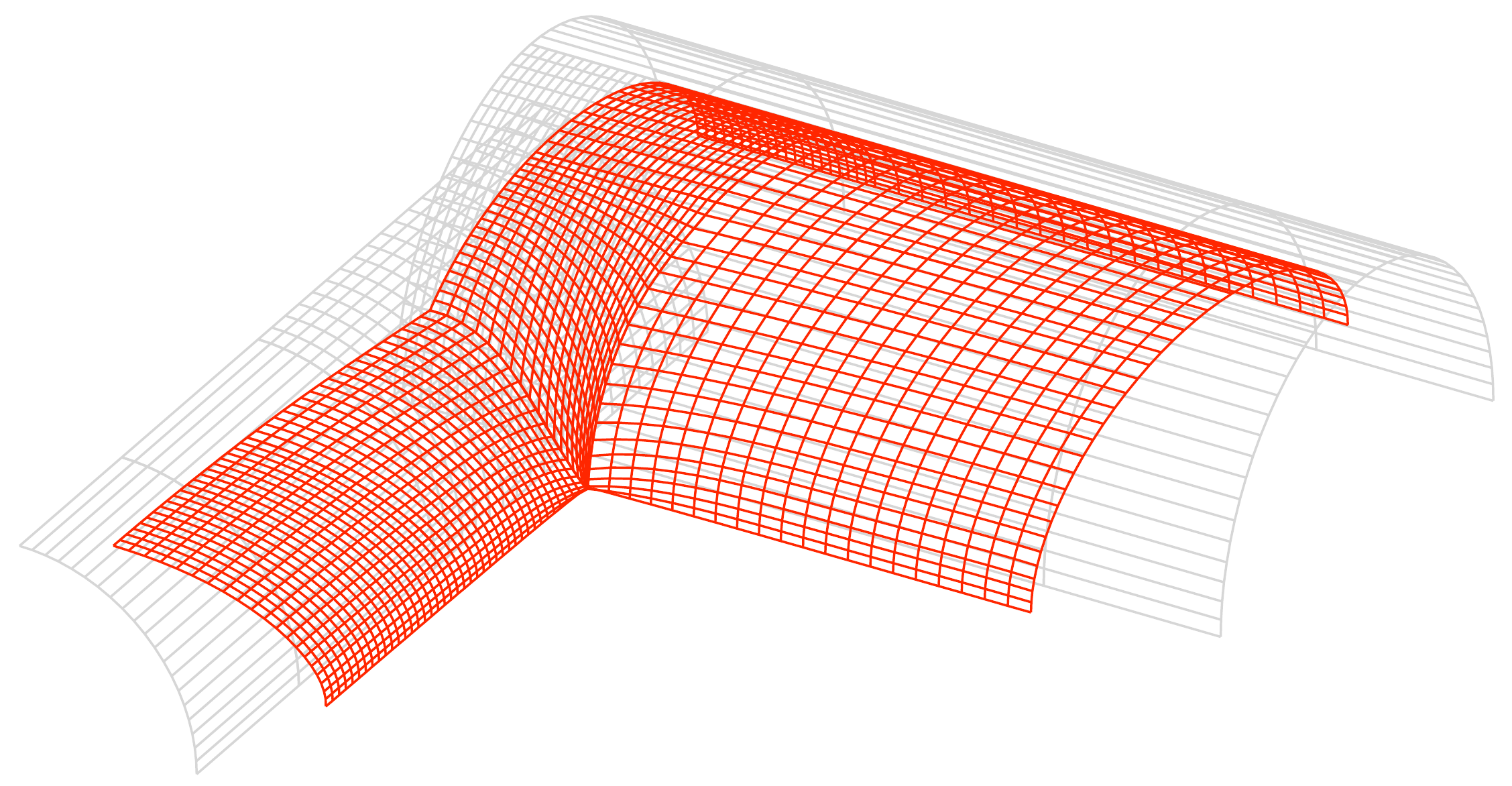}
\caption{Deformed shape of the tunnel}
\label{TB3}
\end{center}
\end{figure}

Using geometry independent field approximation, the basis functions for the definition of the unknown were refined until convergence was achieved. 
The displaced shape of the tunnel walls is depicted in Figure \ref{TB3}.

\begin{figure}
\begin{center}
\includegraphics[scale=0.45]{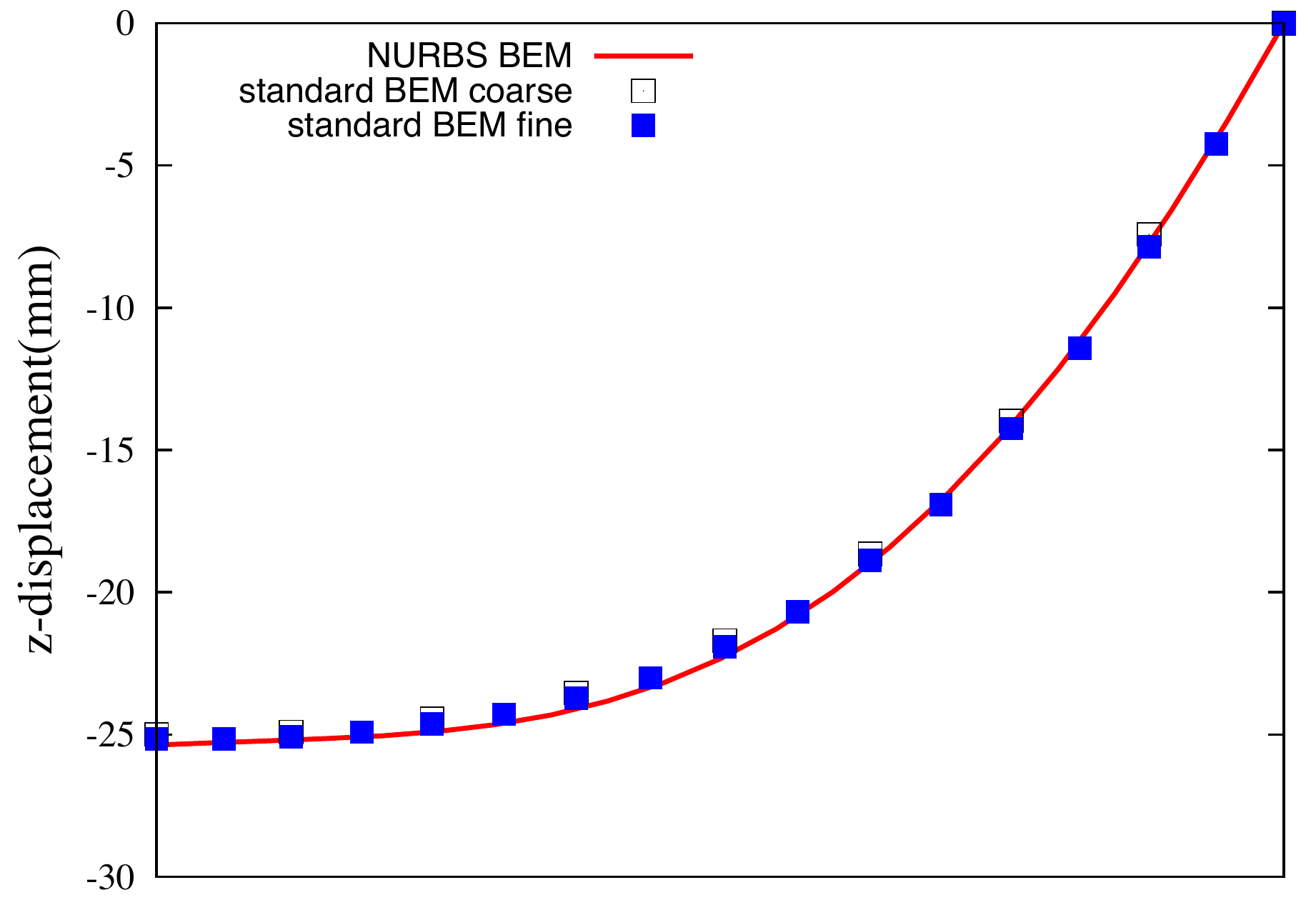}
\caption{Comparison of vertical displacements along the intersection line}
\label{TB4}
\end{center}
\end{figure}

A comparison with a standard isoparametric BEM analysis is shown in Figure \ref{TB4} and good agreement can be observed.

\newpage

\section{Summary and Conclusions}

\begin{figure}
\begin{center}
\includegraphics[scale=1.0]{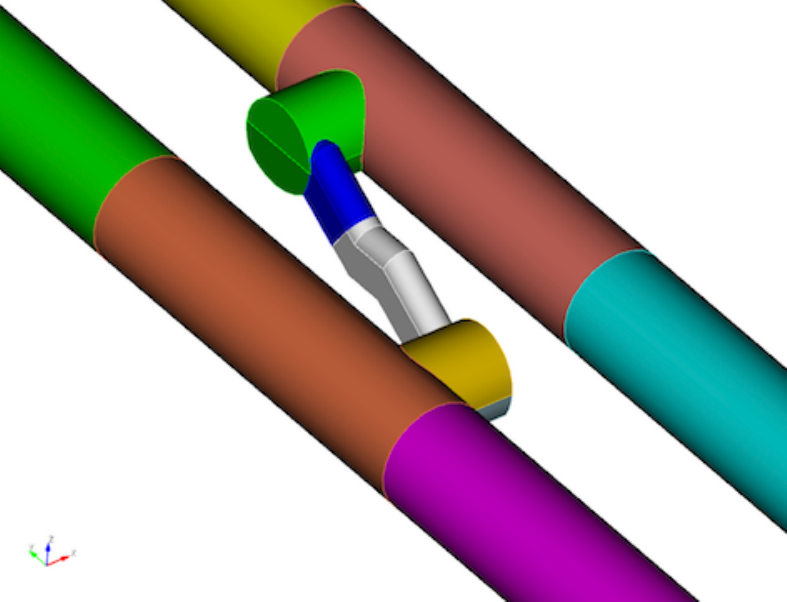}
\caption{Example of a more complex tunneling situation}
\label{Hud}
\end{center}
\end{figure}

A novel approach to trimming was presented. The main advantage of the method lies in its simplicity and ease of implementation into existing software. The method is currently restricted to cases with two trimming curves. However, as has been demonstrated, practical applications exist, where this is the case. Our main interest lies in the application of the BEM to the simulation of underground excavations, where benefits of the method can be seen. 
The geometry of the underground excavation depicted in Figure \ref{Hud} for example is ideally suited to the proposed algorithm. Using a conventional BEM simulation, it has been found difficult to generate a mesh from the CAD data, as elements with bad aspect ratios appeared and had to be repaired. Since the trimming information provided by the CAD program only involves two trimming curves that are not straight lines  in this case, this application lends itself well to the proposed simulation approach.

Extensions of the method are possible. In fact, Gordon-Coons patches (see \cite{Gordon}) provide a framework where four curves can be used to define a surface. One may thus construct an arbitrary patch with four trimming curves in the $u,v$-parameter space. Multiple holes may be considered by using a subdivision of the mapping area.

Trimming plays an important role in realizing the dream of a direct connection between CAD and simulation, without the intermittent step of mesh generation. It is hoped that our contribution provides impetus for a - much needed - increase in research activity in this area.

\section{Acknowledgements}

The work was partially funded by the Austrian Science Fund (FWF) under the project "Fast isogeometric BEM". Thanks are due to T-P. Fries for his discussions and critical comments and to R. Fleissner for making available the MATLAB program for shell analysis.
We also appreciate the suggestions of T-P. Fries and one of the reviewers that Coon's patch technology could be used to extend the capabilities of our method.





\bibliographystyle{elsarticle-num}
\bibliography{ifbbib}







\end{document}